\documentclass{amsart}
\title[A negative answer to Nevanlinna's type question]%
{A negative answer to Nevanlinna'S Type question\\and
a parabolic surface with a lot of\\negative curvature}
\author{Itai Benjamini\and Sergei Merenkov\and Oded Schramm}
\dedicatory{in memory of Bob~Brooks}

\usepackage{amsmath, amssymb, latexsym}
\usepackage{graphicx}
  \ifx\LabelFigloaded\MYundefined\relax
  \else
   \fi

  \def\LabelFigloaded{\relax}

  \chardef\LabelFigCatAt\the\catcode`\@
  \catcode`\@=11

 \let\LabelFigwlog@ld\wlog
 \def\wlog#1{\relax}

 \ifx\\\MYundefined@
    \let\\\relax
 \fi

  \def\ms@g{\immediate\write16}

 \def\N@wif{\csname newif\endcsname }
 \def\Temp@ {\N@wif\ifIN@}
 \ifx\INN@\MYundefined@
    \else \let\Temp@\relax
 \fi
 \Temp@

  \def\IN@{\expandafter\INN@\expandafter}
  \long\def\INN@0#1@#2@{\long\def\NI@##1#1##2##3\ENDNI@
    {\ifx\m@rker##2\IN@false\else\IN@true\fi}%
     \expandafter\NI@#2@@#1\m@rker\ENDNI@}
  \def\m@rker{\m@@rker}
 
  \newtoks\Initialtoks@  \newtoks\Terminaltoks@
  \def\SPLIT@{\expandafter\SPLITT@\expandafter}
  \def\SPLITT@0#1@#2@{\def\TTILPS@##1#1##2@{%
     \Initialtoks@{##1}\Terminaltoks@{##2}}\expandafter\TTILPS@#2@}

 \def\Shifted@@#1#2#3{\setbox0=\hbox{#3}%
   \raise -\dp0\vbox {\kern-#2%
       \hbox {\kern#1\unhbox0\kern-#1}%
           \kern#2}}

 \newcount\gridcount
 \newbox\auxGridbox@ \newbox\hGridbox@ \newbox\vGridbox@
 \newbox\Labelbox@ \newbox\auxLabelbox@
 \newbox\Coordinatebox@
 \newtoks\Labeltoks@
 \newdimen\Wdd@ \newdimen\Htt@
 \newdimen\Wddd@ \newdimen\Httt@
 
 \def\Wr@{\immediate\write16}

 \newdimen\GL@wd
 \def\GridLineWidth#1{\GL@wd=#1}

 \def\gobble#1{}
 \def\EdgeErr@{\Wr@{}%
      \Wr@{\string\Edges\space argument
      1, 10, 100 or 1000 please\string!}%
      }

 \newcount\Edgect@

 \def\Sweepup#1\endSweepup{}

 \def\SetEdges@{%
    \edef\Zr@@s{\expandafter\gobble\number\Edgect@\empty}%
        \count255=0\Zr@@s\relax
        \ifnum\count255=\z@\else\EdgeErr@\show\tailtest\fi
        \count255=1\Zr@@s\relax
        \ifnum\count255=\Edgect@\relax\else\EdgeErr@\show\leadtest\fi
    \EdgGl@b\edef\Zr@s{\expandafter\gobble\Zr@@s\empty}
    \ifnum\Edgect@>\@ne\relax\EdgGl@b\let\L@Dc\empty
        \else\EdgGl@b\edef\L@Dc{\string.}\fi
    \ifnum\Edgect@>\@ne\relax
        \EdgGl@b\edef\Edgescale@##1{\divide##1 by \Edgect@}%
        \else\EdgGl@b\edef\Edgescale@##1{}\fi
    }

 \def\Edges#1{\Edgect@=#1\relax
     \let\EdgGl@b\global \SetEdges@}

 \Edges{1}

 \def\hhrule{\hrule height \GL@wd\vskip-.\GL@wd}

 \def\hRule@{%
   \advance\gridcount -2%
   \vfil\hhrule\vfil
   \llap{\smash{\raise -2.5pt
     \hbox{\L@Dc\number\gridcount\Zr@s\kern2pt}}}%
   \hhrule
   }

\def\vvrule{\vrule width \GL@wd \kern-\GL@wd}

 \def\vRule@{\advance\gridcount 2%
   \hfil\vvrule\hfil
   \setbox\auxGridbox@=\vbox to 0pt
      {\vskip \Htt@\vskip 2pt
        \hbox to 0pt{\hss\L@Dc\number\gridcount\Zr@s\hss}\vss}%
      \wd\auxGridbox@=0pt \box\auxGridbox@
   \vvrule
   }

 \def\PlaceGrid@@{\gridcount=10 
  \setbox\hGridbox@=\hbox{%
        \hbox{%
             \hskip-.4pt\vrule
             \vbox to \Htt@{%
               \offinterlineskip\parindent=\z@\relax
               \hbox to \Wdd@{\hfil}
               \hRule@\hRule@\hRule@\hRule@
               \vfil\hhrule\vfil}%
             \vrule\hskip-.4pt}
    }%
  \gridcount=0%
  \setbox\vGridbox@=\hbox{%
      \vbox{\offinterlineskip\parindent=0pt\hsize=0pt
         \vskip-.4pt\hrule%
         \hbox to \Wdd@{%
                 \vtop to \Htt@{\vfil}%
                 \vRule@\vRule@\vRule@\vRule@
                 \hfil\vvrule\hfil}%
         \hrule\vskip-.4pt}}%
  \wd\hGridbox@=0pt\ht\hGridbox@=0pt
  \wd\vGridbox@=0pt\ht\vGridbox@=0pt
  \hbox{\box\hGridbox@\box\vGridbox@}%
  }

 \def\LabelsGlobal{\def\LabGl@b{\global}}
 \def\LabelsLocal{\def\LabGl@b{}}
 \LabelsGlobal 

 \def\SetLabels#1\endSetLabels{%
   \LabGl@b\Labeltoks@={#1()\\}%
   }

 \LabGl@b\Labeltoks@={()\\}

 \def\ShowGrid{\LabGl@b\let\PlaceGrid@\PlaceGrid@@}
 \def\HideGrid{\LabGl@b\let\PlaceGrid@\relax}
 \def\Grids{\ShowGrid\LabGl@b\let\GridSwitch@\ShowGrid}
 \def\noGrids{\HideGrid\LabGl@b\let\GridSwitch@\HideGrid}

 \noGrids

 \def\bAdjust@@{%
     \setbox\auxLabelbox@=\hbox{\raise \dp\auxLabelbox@
            \box\auxLabelbox@}}
 \def\bAdjust@{\let\vAdjust@\bAdjust@@}

 \def\eAdjust@@{\dimen0=-.5\ht\auxLabelbox@
     \advance\dimen0 by .5\dp\auxLabelbox@
     \setbox\auxLabelbox@=
            \hbox{\raise\dimen0\box\auxLabelbox@}}
 \def\eAdjust@{\let\vAdjust@\eAdjust@@}

 \def\tAdjust@@{%
     \setbox\auxLabelbox@=\hbox{\raise-\ht\auxLabelbox@
            \box\auxLabelbox@}}
 \def\tAdjust@{\let\vAdjust@\tAdjust@@}

 \let\vAdjust@\relax

 \def\lAdjust@{\let\hAdjust@\rlap}
 \def\rAdjust@{\let\hAdjust@\llap}

 \let\hAdjust@\relax\let\vAdjust@\relax

 \def\FetchLabel@#1(#2)#3\\{%
     \IN@0#2@@\ifIN@
        \setbox0=\hbox{\ignorespaces#1#3\unskip}%
        \ifdim\wd0>0pt
           \ms@g{}%
           \ms@g{ !!! Bad label(s)? !!!}%
           \message{ #1(#2)#3}%
        \fi
        \def\LabelMole@##1\endFetchLabel@{%
            \IN@0()\\@##1@%
            \ifIN@\def\Temp@{\FetchLabel@##1\endFetchLabel@}%
            \else\def\Temp@{}%
            \fi
            \Temp@
           }%
     \else
       \ignorespaces#1\unskip
       \setbox\auxLabelbox@=%
         \hbox to 0pt{\hss\ignorespaces\hAdjust@
          {\ignorespaces#3\unskip}\hss}%
       \vAdjust@
       \let\hAdjust@\relax\let\vAdjust@\relax
       \AugmentLabelBox@@{#2}%
       \ht\Labelbox@=0pt\dp\Labelbox@=0pt
       \let\LabelMole@\FetchLabel@%
     \fi\LabelMole@}

 \newtoks\XYSep@ 
 \def\SetXYSeparator#1{%
     \IN@0#1@@\ifIN@\XYSep@{*}%
     \else
     \XYSep@{#1}%
     \fi
     }

 \SetXYSeparator*

 \def\AugmentLabelBox@@#1{%
     \IN@0\the\XYSep@ @#1@\ifIN@
       \SPLIT@0\the\XYSep@ @#1@%
       \setbox\Labelbox@=\hbox to 0pt{%
         \unhbox\Labelbox@
         \Shifted@@{\the\Initialtoks@\Wddd@}%
         {\the\Terminaltoks@\Httt@}%
         {\box\auxLabelbox@}}%
     \else
         \ms@g{}%
         \ms@g{ !!! Bad insertion point. !!!}%
         \message{ (#1\ this point was rejected.)}%
     \fi
    }

 \def\FetchOption@#1[#2]#3\endFetchOption@{%
    \def\temp{#1}
    \ifx\temp\empty
       \Edgect@=#2\relax
       \let\EdgGl@b\relax
       \SetEdges@
       \Cleaner@#3%
    \fi}

 \def\Cleaner@#1[@]{\Labeltoks@{#1}}
     
 \def\PlaceLabels@@{\mathsurround=0pt
     \def\Cr@{\\}%
     \let\L\lAdjust@\let\R\rAdjust@
     \let\B\bAdjust@\let\E\eAdjust@\let\T\tAdjust@
     \expandafter\FetchOption@\the\Labeltoks@[@]\endFetchOption@
     \Wddd@=\Wdd@ \Edgescale@\Wddd@ 
     \Httt@=\Htt@ \Edgescale@\Httt@
     \expandafter\FetchLabel@\the\Labeltoks@\endFetchLabel@
     \box\Labelbox@
     }%

 \let \PlaceLabels@\PlaceLabels@@

 \def\AffixLabels#1{\setbox\Coordinatebox@=\hbox{#1}%
      \Wdd@=\wd\Coordinatebox@ \Htt@=\ht\Coordinatebox@
      \advance\Htt@ \dp\Coordinatebox@
      \hbox{\copy\Coordinatebox@\kern-\Wdd@ 
           \Shifted@@{0pt}{-\dp\Coordinatebox@}%
           {\PlaceLabels@\PlaceGrid@}%
           \kern\Wdd@}%
      \GridSwitch@ 
      \LabGl@b\Labeltoks@{()\\}%
      }
 
   \let\wlog\LabelFigwlog@ld   
   \catcode`\@=\LabelFigCatAt  

                                By

              Raymond S\'eroul <A18645@FRCCSC21.BITNET>
                                and 
              Laurent Siebenmann <lcs@topo.math.u-psud.fr>
    
              VERSIONS: July 1991, Oct 1991, Jan 1992, July 1992

INTRODUCTION

      This labelling package is intended for TeX users who
rely on non-TeX sources for for their graphics inserts.  It
provides means for adding TeX labels to such inserts with a
minimum of fuss. 

       For most labels, TeX users have in the past found it
reasonably convenient to rely on non-TeX sources. Typical
occasions when an inescapable need for TeX labels seemed to
arise are

 (a) when the graphics program lacks certain exotic or complex
mathematical symbols

 (b) when the very highest typographical quality is wanted for the
labels

 (c) when labels included with the graphics fail to print, 
 and you cannot figure out why (cf. boxedeps.doc).  The labels
 provided by labelfig.tex are 100

       Since this package first appeared, many users, who in the
past scarcely dreamed of using TeX labels, have come to use
nothing but.  So it is now appropriate to add

Intoxication Warning:  TeX labels may be addictive and expensive. 

     If you have a fast preview you may disagree, and even find
that this package provides an agreeable paste-up environment; see
extra applications at end.

     Note to publishers: It is possible and convenient to ultimately
export the TeX labels produced by labelfig.tex to become an integral
part of the EPS file. This is often desired by a publisher who typically
uses an "upmarket" graphics or page layout program, with which the
staff is skilled in perfecting figures.  See Appendix I for
a recipe.

     The authors are grateful to Patrick Ion of Math Reviews for
helpful comments and encouragement.

BASIC INSTRUCTIONS

    After reading in the macro file using

preview or proof your figure with a coordinate grid printed on
top, by typing the following:

    \ShowGrid  
    \AffixLabels{<the graphics insertion>}

Here <the graphics insertion> is what you would type to insert
the graphics object alone without the grid.  This must provide
for the space around it. For example <the graphics insertion>
might well be \BoxedEPSF{MyFigure scaled 700} using the
boxedeps.tex macro package (from same source); this provides a
TeX box containing the encapsulated PostScript insert specified by
the file MyFigure. \AffixLabels{...} provides the grid (supposing
\ShowGrid is present) and later, once you have specified labels
using the grid, it will "tack on" the labels.

     The grid is a sort of (usually elongated) checkerboard of
ten rows and ten columns and its (internal) partitions are by
default numbered  .1, ... ,.9  both horizontally (X-coordinate
running left to right) and vertically (Y-coordinate running bottom
to top).  Thus the points enclosed by the grid correspond to the
points of the unit square in the cartesian "X-Y" plane, the lower
left corner corresponding to the origin (0,0).  By extrapolation,
the full page corresponds to a larger rectangle in the plane.

     These coordinates serve to position labels as follows.
Before the \AffixLabels{...} command type label specifications:

  \SetLabels
   (<X-coordinate>*<Y-coordinate>) <first label> \\
   .
   .
   .
   (<X-coordinate>*<Y-coordinate>)  <last label> \\
  \endSetLabels

Each row specifies one label and is terminated by \\.  In each
row, the position indicator comes first; it is written as a
standard cartesian point except that the X- and Y- coordinates
are separated by * rather than a comma because TeX allows a
comma as decimal point. There are no dimension units to specify
as the unit is the grid itself.

     By default, this cartesian point specifies where the middle
of the baseline of the label will be located.  However if you precede
the point by \L [or \R] the left [or right] edge of the baseline will
be located there. Similarly you may also precede the point by \T, \E,
or \B to vertically align the top equator or bottom of the label box
at the specified point.  This gives nine standard positions of
the label with respect to the insertion point --- corresponding to
the eight principle points of the compas and the center

                     \L\T     \T      \R\T

                     \L\E     \E      \R\E

                     \L\B     \B      \R\B

But this neglects the default "baseline" level of TeX,
giving potentially three more positions

                     \L    <no tag>   \R

For text, the baseline level is often the preferred. Its relation to
the others is variable. It will often coincide with the bottom level,
as happens for "X".  But it is often distinct, as for "g", in which
case you have in all 12 distinct positions rather than 9.

     It is convenient to think of this specification of label
position as attaching the label by a thumb-tack to the coordinate
grid. There are up to twelve positions of the thumb-tack on the
label, while the position of the thumb-tack on the coordinate grid is
arbitrary.  Normally, one choses the position of the thumb-tack on
the label to be the one that is the closest to the item being
labeled.  There are good reasons for this "rule of thumb":

   (a)  It facilitates correct positioning at first try.

   (b)  If the scale of the figure must be altered after labels
have been affixed, the labels have a good chance of remaining well
positioned.

   (c)  The visible grid need not extend beyond the "bounding box"
for the figure, because the best preferred position is always
(at least almost) within the bounding box .

The second reason is particularly important. Indeed it often
happens that scale has to be altered after labelling begins, in
order to either provide space for the labels, or to adjust
proportions between the labels and the figure.  (The size of labels
is unaffected by scaling.)

     Here is an artificial but self-contained test which uses
TeX rules to make a graphics object.

TEST

    Do not skip this!

 \def\FrameIt#1{\hbox{\vrule$\vcenter {\hrule\kern3pt%
             \hbox {\kern3pt #1\kern3pt}%
               \kern3pt\hrule}$\relax\vrule}}

 \def\Caption#1#2{\FrameIt{%
       \vtop {\hsize=#1\relax \parindent=0pt
         \leftskip=0pt \rightskip=0pt plus15pt
         \parfillskip=0pt
         \lineskip=1pt\baselineskip=0pt
         #2}}}

 \def\FirstQuadrant{\hbox to 100pt{\vrule\vbox to 100pt{%
        \hbox to 100pt{\hfil}\vfil\hrule}\hss}}

  \SetLabels
    \R(.5*.2) $\zeta\,\cdot$\\
    (.9*-.10) $\xi$\\
    \R(-.03*.9) $\eta$\\
    \T(.5*.9) \Caption{70pt}{%
          \it The norm of
          $g(\xi+i\eta)$ is indicated on
          contours of this invisible surface.}\\
  \endSetLabels

  \AffixLabels{\FirstQuadrant} \end

  Note that the coordinates to use for labels are indicated on the
edges of the grid (when visible) corresponding to the conventional
x- and y- axes of the Cartesian plane. By default the grid is
1-by-1. However, by the command \Edges{100}, you can change this
to 100-by-100 and many users find this alternative most
convenient. Place the command \Edges{...} in your style file (or
header) since its effect is is global. Other possible edge values
are 10 and 1000.

  If you use the command \Edges{...} at all, do so with care.  For
if you accidentally delete an \Edges{...} command your labels will
abruptly be badly misplaced and may logically but mysteriously
generate "dimension too big" errors under TeX and "off page" errors
under your driver.  

  You can dictate the edgescale for an individual figure by giving
the scale in brackets immediately after \SetLabels.  Thus, to
import into an article using say \Edge{100} a figure labelled using
another edgescale, say the original 1-by-1 default, you can use
\SetLabels[1]...\endSetLabels.

GETTING IT DOWN PAT

     Complicated labeling deserves the same respect as
complicated mathematics.  Do not expect it to come out perfect the
first time!  What is needed in either case is a mechanism to
repeatedly typeset troublesome pieces.

     One mechanism is always available.  One does complicated
labelling in a separate "test" file involving just the figure being
labelled;  a texpert will know how to \dump TeX's current state as
a temporary format that restarts rapidly at each retry.  Usually,
one then pastes the completed labelled figure back into the main
TeX file, but, of course, one can also \input it as an auxiliary
file.

     If you do not have a TeXpert at handy, here is a first
approximation to an efficient setup. By deletions reduce a copy
of your article to just a few lines before and after the figure.
Now label the figure, and finally, copy and paste the labelled
figure to the original article. Then copy the next figure to label
into this testbed and repeat. The TeXpert can improve the  speed
at which TeX starts up, by compiling a format specifically for
your article; just one caution: best NOT include in the format
ephemeral details of setup like \Set<mydriver>ArtSpecials (from
boxedeps.tex because this reads  figure dimensions which you may
change during your work session.

     An improved mechanism to repeatedly typeset troublesome
pieces is now available on the Macintosh; it is called LinoTeX;
see the same ftp sources.  It could be set up on many types
of computer.

     Before using labelfig.tex to attach labels to a graphics
object inserted using boxedeps.tex or BoxedArt.tex, make it a
firm rule to carefully adjust the bounding box using the trimming
commands of these packages, and also at least tentatively scale
and position the object. Beware of changing the grid inadvertently
after the labels have been positioned.  For example, correcting
the bounding box of a PostScript graphics object can foul up the
labels by changing the coordinate grid to which the labels are
attached. This is particularly true for the trimming  commands of
boxedeps.tex and BoxedArt.tex. However, as noted already, change
of scale is much less disruptive, and modest adjustments should be
well tolerated.

     Sometimes the labels protrude so far from the bounding box
of a figure that the figure has to be repositioned.  Best do this
by ad hoc spacing, say using \hglue and \vglue; altering the
bounding box would create a vicious circle.

     Remember that you are responsible for preventing labels
from overlapping. You are responsible for all label typography
including size and style. A label is really just about anything
that can be put in a TeX box. Note that spaces at the beginning
and end of labels will normally be suppressed; if you really want
them you must protect them with TeX braces.

     This package temporarily sets the \mathsurround parameter
of TeX to zero  while the labels are being affixed. This is done
because nonzero \mathsurround space would influence the position
of left and right aligned labels; then, when a texpert or printer
modifies mathsurround, diagram labeling might be disastrously
altered. There is a small price to pay involving labels that are
formatted as caption boxes including mathematics: you  may want or
need to specify an explicit mathsurround space within the caption
box; it will not influence anything outside.

     Those hostile to the use of * as separator between
the X and Y coordinates of label insertion points, are free to
impose another using \SetXYSeparator{<the new separator>}.  
Americans may prefer "," to "*" since they never use a 
comma as a decimal point; on the other hand, * may be more visible.

APPENDIX (I)  MERGING labelfig.tex LABELS INTO AN EPSF GRAPHICS OBJECT.

     As promised in the introduction, here is a recipe useful for
publishers. It works at least on Macintosh and at least for vectorized
graphics and Adobe type1 fonts.  (There is surely a similar recipe for
PCs under MSWindows.)

 (a)  Use boxedeps.tex utility to integrate the figure given by the eps
file, "x.eps" say, with a visible frame around it.  See
\ShowDisplacementBoxes command in boxedeps.tex.  To get precise results
automatically it is important to use the \Trim... commands of
boxedeps.tex making the "DisplacementBox" neatly fit the figure.

 (b)  Use the TeX printer driver and LaserWriter (versions >= 8.1.1) to
export to an EPSF the DVI page containing the integrated, labelled
figure. You now have an EPS file  "xx.eps"  that contains too much, and at
the wrong scale, and at wrong position.

 (c)  Convert the EPSF to an Adode Illustrator format EPSF using
the shareware utility called epsConvert by Sam Weiss
1993-- (currently $25).

 (d)  In Illustrator (or a compatible program), group the labels and the
"DisplacementBox"; copy them to the clipboard and paste them into "x.ps".
This step requires that all the label fonts be "visible to the Macintosh.

 (e)  Translate and scale the pasted group consisting of the labels plus
the "DisplacementBox" so as to make the "DisplacementBox" the bounding
box of (labelless) figure represented by "x.eps".  At this point the
labels will be correctly placed on the figure "x.eps".

 (f)  Ungroup and delete the "DisplacementBox".  The result is the
desired single EPS file, "x+.eps" say, It contains the original figure
plus its labels.  

     Using grouping and ungrouping appropriately in "x+.eps", a
publisher's staff can very efficiently improve label positions etc.

APPENDIX II)  SOME EXOTIC APPLICATIONS

     The grid of labelfig.tex is analogous to a light-table in
classical page makeup with wax or latex glue.  In principle, you
can use it to compose any page from its indivisible parts.  This
even has some of the artisanal charm of classical paste-up
provided you have a fast screen preview to make the process
"interactive".

     In practice labelfig.tex is a tool for nonstandard jobs.
Here are a few going beyond the labelling already discussed.

(I)  GRAPHICS INTEGRATION.

     This is accomplished by treating the imported graphics
objects as labels.  The underlying graphics object is then
typically an empty  \vbox to <dimension>{\vfill} in a TeX
\midinsert...\endinsert construction.  A label line
might be of the form

   (.1*.1) \\

The exact form of the special command varies from driver to
driver.  However, in the case of encapsulated PostScript graphics
(EPSF norm), by relying on boxedeps.tex, one can have the
following standard syntax (independant of driver  (see
boxedeps.doc for details.
  
  (.1*.1) \BoxedEPSF{MyFigure scaled <scale in mils>}\\

This may be slow since it requires TeX to read the PostScript
file to read bounding box using many complex macros.  So you
may want to try

  (.1*.1) \EPSFSpecial{MyFigure}{<scale in mils>}\\

which is fast and driver independant, but it squashes the
bounding box, normally to its lower left corner.

     Similarly for graphics of the Macintosh PICT norm ---
using BoxedArt.tex (same sources) in place of boxedeps.tex.

     This approach to integration is to be recommended when
one is assembling a composite graphics object.

 (II)  COMMUTATIVE DIAGRAM ENHANCEMENT

     Commutative diagrams or arrays of mathematical objects
connected by arrows of various sorts are common in mathematics.
The mathematical objects require the use of TeX.  Recently TeX
acquired a good collection of arrows of all slopes --- that of
LamSTeX --- plus pwerful macros to build the diagrams.

     However, even the LamSTeX collection is often
inadequate; it lacks for example double shafted arrows, dotted
arrows and curved arrows. Fortunately it is possible to produce
such arrows on an individual basis using sophisticated graphics
programs such as Illustrator and AldusFreehand (both serving
the EPSF norm) or using Metafont (with its public domain norm).
Since the creation of each new arrow is a work of love, you
probably want to limit the number of arrows by using LamSTeX
for most arrows. The 40K commutative diagram module of LamSTeX
has been adapted to work with AmSTeX and a copy may be posted
with LabelFig and related files. Unfortunately no one has yet
offered a version that works with Plain TeX or LaTeX.

       Suffice it here to say that when the exotic arrow has
been somehow imported into TeX, labelfig.tex treats it as a
label that one affixes to the commutative diagram.  Two other
steps will be treated in separate notes, namely the matter of
extracting the dimension specifications for the arrow and the
construction of the arrow --- for these steps are far from
unique and often depend intimately on your computer environment. 
Notes for the Macintosh-Textures-Illustrator combination are
found in the file ExoticArrows.doc.

 (III) NESTING 

Ingenuity pays off in exploiting labelfig.tex. One can
mix graphics and typography quite freely.  labelfig.tex is good
for freeform or overlapping arrangements, while boxedeps.tex (or
BoxedArt.tex) is best for regimented non-overlapping
arrangements --- and the two can be combined.

     The default behavior of labelfig.tex is not ideal 
for nesting objects, because to prevent trouble for beginners
the register for labels is globally cleared when \AffixLabels
concludes.  But there are switches available

      \LabelsGlobal      \LabelsLocal

which change this.  To understand this, extend the above test 
by something like:

 \LabelsLocal

 \SetLabels
    (.5*.5) AAA\\
 \endSetLabels

 {
 \SetLabels
    (.5*.5) ZZZ\\
 \endSetLabels
   \AffixLabels{\FirstQuadrant}
 }

   \AffixLabels{\FirstQuadrant}

     There are however potential pitfalls.  Neither
labelfig.tex nor boxedeps.tex has been tested under extreme
conditions. Problems may occur if their procedures are
indiscriminately nested. For boxedeps.tex (not labelfig.tex)
there is a precise cause for worry, namely many of its
variables are "global", which means that TeX braces will not
provide the protection one might expect.

COMMAND SUMMARY FOR labelfig.tex

  Here [...] means optional (one or zero)
       [...]* means any number of such constructs

  \SetLabels
    [[<P>](<X><Sep><Y>) <label> \\]*
  \endSetLabels
  \ShowGrid  
  \AffixLabels{<the figure>}

   --- <P> is tack position, one of eleven or empty
              order irrelevant

                   \L\T      \T      \R\T

                   \L\E      \E      \R\E

                     \L               \R

                   \L\B      \B      \R\B

   --- (<X><Sep><Y>) insertion point;
  <Sep> is separator, = * by default;
  \SetXYSeparator{<Sep>} changes it.
   <X> and <Y> are real numbers

  --- <label> a label to attach 

  --- <the figure> the figure to label 

  \GlobalLabels (default)     
  \LocalLabels  setting for nested constructs.

 \Grids makes ALL grids appear; \HideGrid then makes just next disappear.
 \noGrids returns to default.  The commands are always global.

 \GridLineWidth{<dimension>} adjusts width of grid lines. Default is very
small, to give "hairline" effect. If your grid lines are missing try
setting \GridLineWidth{1pt}.

 \Edges#1 globally changes the edge size of all grids to the numerical 
value #1, which must be 1, 10, 100, or 1000.  The default is 1.

VERSION HISTORY.
 --- Jan 1993: \Edges#1 and [??] option after \SetLabels
 --- July 1992: \Grids, \noGrids, \HideGrid;
       Gridlines become hairlines; \GridLineWidth{<dimension>}.
 --- Oct 1991, Jan 1992: \SetXYSeparator{<Sep>},  \LabelsGlobal,
       \LabelsLocal.
 --- July 1991: first release

Address for bugs and other feedback:

        Raymond S\'eroul
        IREM and Lab. de Typographie Informatise
        Univ. Rene Descartes
        Strasbourg

    Tel 33-88-41-63-45
    Email:  A18645@FRCCSC21.BITNET

        Laurent Siebenmann
        Mathematique, Bat. 425,
        Univ de Paris-Sud,
        91405-Orsay,
        France

    Tel 33-1-6941-7949; 
    Email: lcs@topo.math.u-psud.fr

\newcommand\C{{\mathbb C}}

\newcommand\D{{\mathbb D}}

\newcommand\R{{\mathbb R}}

\def\area{\mathop{\mathrm{area}}}
\def\length{\mathop{\mathrm{length}}}
\newcommand\OC{{\overline{\mathbb C}}}

\newtheorem*{definition}{Definition}

\begin{document}

\thanks{Research of the second author supported by NSF grant DMS-0072197}
\address{Department of Mathematics\\ Weizmann Institute of Science\\
Rehovot 76100, Israel} \email{itai@math.weizmann.ac.il}
\address{Department of Mathematics\\ Purdue University\\ West Lafayette,
Indiana 47907} \email{smerenko@math.purdue.edu}
\address{Microsoft Research\\ One Microsoft Way\\ Redmond, WA 98052}
 \email{schramm@microsoft.com}

\abstract{
Consider a simply connected Riemann surface represented by a Speiser graph.
Nevanlinna asked if the type of the surface is determined by the mean
excess of the graph:
whether mean excess zero implies that the surface is parabolic and
negative mean excess implies that the surface is hyperbolic.
Teichm\"uller gave an example of a hyperbolic
simply connected Riemann surface whose mean excess is
zero, disproving the first of these implications.
We give an example of a simply connected parabolic Riemann surface with
negative mean excess, thus disproving the other part.
We also construct an example of a complete,
simply connected, parabolic surface with nowhere positive
curvature such that the integral of curvature in any disk
about a fixed basepoint is less than $-\epsilon$
times the area of disk, where $\epsilon>0$ is some constant.}
\endabstract

\maketitle

\section{Introduction}\label{S:Intro}

The uniformization theorem states that for every simply-connected
Riemann surface $X$ there exists a conformal map $\varphi:\ X_0\to
X$, where $X_0$ is either the complex plane $\C$,
the open unit disc $\D=\{z\in\C:|z|<1\}$, or the
extended complex plane (Riemann sphere) $\OC=\C\cup\{\infty\}$,
and these possibilities are mutually exclusive~\cite{lA73}.
The map $\varphi$ is called the {\emph{uniformizing
map}}. A simply-connected Riemann surface $X$ is said to have
{\emph{hyperbolic}}, {\emph{parabolic}}, or {\emph{elliptic}}
{\emph{type}}, according to whether it is conformally equivalent
to $\D$, $\C$, or $\OC$, respectively. In what follows, we assume
that $X$ is not compact, thus excluding the elliptic case from
consideration.

We are interested in the application of the Uniformization Theorem
to the following construction. A {\emph{surface spread over the
sphere}} is a pair $(X, \psi)$, where $X$ is a topological surface
and $\psi:\ X\to\OC$ a continuous, open and discrete map. The map $\psi$
is called a {\emph{projection}}. Two such surfaces $(X_1, \psi_1),\
(X_2, \psi_2)$ are {\emph{equivalent}}, if there exists a
homeomorphism $\phi:\ X_1\to X_2$, such that $\psi_1= \psi_2\circ\phi$.
According to a theorem of Sto\"{\i}low \cite{sS56}, a continuous
open and discrete map $\psi$ near each point $z_0$ is homeomorphically
equivalent to a map $z\mapsto z^k$. The number $k=k(z_0)$ is called
the {\emph{local degree}} of $\psi$ at $z_0$. If $k\neq 1$, $z_0$ is
called a {\emph{critical point}} and $\psi(z_0)$ a {\emph{critical
value}}. The set of critical points is a discrete subset of $X$.
The theorem of Sto\"{\i}low implies that there exists a unique
conformal structure on $X$ which makes $\psi$ into a meromorphic
function. If $X$ is simply-connected, what is the type of the
Riemann surface so obtained? This is one version of the type
problem. Equivalent surfaces have the same type.

Rolf Nevanlinna's problem concerns a particular class of surfaces
spread over the sphere, denoted by $F_q$. Let $\{a_1, \dots,
a_q\}$ be distinct points in $\OC$.

\begin{definition}\label{D:Fq} A surface $(X, \psi)$
belongs to the class $F_q=F(a_1, \dots, a_q)$, if $\psi$ restricted
to the complement of $\psi^{-1}\bigl(\{a_1,\dots,a_q\}\bigr)$
is a covering map onto its image $\OC\setminus\{a_1,\dots,a_q\}$.
\end{definition}

Assume that $(X,\psi)\in F_q$ and $X$ is noncompact.
We fix a Jordan curve $L$, visiting the points $a_1,\dots, a_q$
in cyclic order.
The curve $L$ is usually called a {\emph{base curve}}. It
decomposes the sphere into two simply-connected regions $H_1$,
the region to the left of $L$, and $H_2$, the region to the right of
$L$.
Let $L_i$, $i=1,2,\dots,q$, be the arc of $L$
from $a_i$ to $a_{i+1}$ (with indices taken
modulo $q$). Let us fix
points $p_1$ in $H_1$ and $p_2$ in $H_2$, and choose $q$
Jordan arcs $\gamma_1,\dots, \gamma_q$ in ${\OC}$, such that each
arc $\gamma_i$ has $p_1$ and $p_2$ as its endpoints, and has a
unique point of intersection with $L$, which is in $L_i$.
We take these arcs to be interiorwise disjoint, that is,
$\gamma_i\cap\gamma_j=\{p_1,p_2\}$ when $i\ne j$.
Let $\Gamma'$ denote the graph embedded in $\OC$, whose vertices are
$p_1,\ p_2$, and whose edges are $\gamma_i,\ i=1,\dots, q$, and let
$\Gamma=\psi^{-1}(\Gamma')$. We identify $\Gamma$ with its image in
$\R^2$ under an orientation preserving homeomorphism of $X$ onto $\R^2$.
The graph $\Gamma$ has the following properties:
\begin{enumerate}
\item $\Gamma$ is infinite, connected;
\item $\Gamma$ is homogeneous of degree $q$;
\item $\Gamma$ is bipartite.
\end{enumerate}
A graph, properly embedded in the plane and satisfying these properties
is called a {\emph{Speiser graph}}, also known as a
{\emph{line complex}}. The vertices of a Speiser graph $\Gamma$
are traditionally marked by $\times$ and $\circ$, such that each edge
of $\Gamma$ connects a vertex marked $\times$ with a vertex marked $\circ$.
Such a marking exists, since $\Gamma$ is bipartite.
Each {\emph{face}} of $\Gamma$, i.e.,  connected component of
$\R^2\setminus\Gamma$, has either a finite even number of edges along its
boundary, in which case it is called an
{\emph{algebraic elementary region}},
or infinitely many edges, in which case it is called a
{\emph{logarithmic elementary
region}}. Two Speiser graphs $\Gamma_1, \ \Gamma_2$ are said to be
{\emph{equivalent}}, if there is a sense-preserving homeomorphism
of the plane, which takes $\Gamma_1$ to $\Gamma_2$.

The above construction is reversible. Suppose that the faces of a
Speiser graph $\Gamma$ are labelled by $a_1, \dots, a_q$, so that
when going counterclockwise around a vertex $\times$, the indices
are encountered in their cyclic order, and around $\circ$ in the
reversed cyclic order.
We fix a simple closed curve $L\subset\OC$ passing through $a_1, \dots, a_q$.
Let $H_1,H_2,L_1,\dots,L_{q}$ be as before.
Let $\Gamma^*$ be the planar dual of $\Gamma$.
If $e$ is an edge of $\Gamma^*$ from a face of $\Gamma$ marked
$a_j$ to a face of $\Gamma$ marked $a_{j+1}$, let $\psi$ map
$e$ homemorphically onto the corresponding arc $L_j$ of $L$.
This defines $\psi$ on the edges and vertices of $\Gamma^*$.
We then extend $\psi$ to the faces of $\Gamma^*$ in the obvious
way.  This defines a
surface spread over the sphere $(\R^2, \psi)\in F(a_1, \dots, a_q)$.
See~\cite{Nevbook} for further details.

For a Speiser graph $\Gamma$, Nevanlinna introduces the
following characteristics. Let $V(\Gamma)$ denote the set of
vertices of the graph $\Gamma$. To each vertex $v\in V(\Gamma)$ we
assign the number
$$
E(v)= 2-\sum_{f\in F(v)} (1-1/k_f),
$$
where $F(v)$ denotes the set of faces containing $v$ on their boundary
and $2\,k_f$ is the number of edges on the boundary of $f$,
$k_f\in\{1,2,\dots,\infty\}$.
The function $E:\ V(\Gamma)\to \R$ is called the {\emph{excess}}.
This definition is motivated by considering the curvature, as follows. The
$\psi$-pullback of the spherical metric $2|dw|/(1+|w|^2)$ is generally
singular, i.e., it may degenerate on $\psi^{-1}(\{a_1,\dots,a_q\})$. The
surface $X$, endowed with the pullback metric, is a
spherical polyhedral surface, which is a particular kind of orbifold.
The {\emph{integral curvature}} $\omega$ on $X$ is a signed Borel
measure, so that for each Borel subset $B\subset X$, $\omega(B)$
is the area of $B$ with respect to the pullback metric minus
$2\pi\sum_z(k_z-1)$, where the sum is over all critical points $z\in
B$ and $k_z$ is the local degree of $\psi$ at $z$.

Each vertex of $\Gamma$ represents a hemisphere, and each face of
$\Gamma$ with $2k$ edges, represents a critical
point, where $k$ is the local degree of $f$ at this point. Therefore,
each vertex of $\Gamma$ has positive integral curvature $2\pi$,
and each face with $2k$ edges has negative integral curvature
$-2\pi(k-1)$. We assign the negative curvature evenly to all the
vertices of the face. A face with infinitely many edges
contributes $-2\pi$ to each vertex on its boundary. The curvature
assigned to every $v\in V(\Gamma)$ is exactly $2\pi E(v)$.

Nevanlinna also defines the mean excess of a Speiser graph
$\Gamma$. We fix a base vertex $v_0\in V(\Gamma)$, and consider an
exhaustion of $\Gamma$ by a sequence of finite graphs
$\Gamma_{(i)}$, where $\Gamma_{(i)}$ is the ball of combinatorial
radius $i$, centered at $v$. By averaging $E$ over all the
vertices of $\Gamma_{(i)}$, and taking the limit, we obtain the
{\emph{mean excess}}, denoted $E_m=E_m(\Gamma)$, provided that
the limit exists.
If the limit does not exist, we consider the {\it{upper mean
excess}} $\overline E_m$ and {\it
{lower mean excess}} $\underline E_m$, which are
the upper ($\limsup$) and lower ($\liminf$) limits, respectively.

\medskip\noindent{\bf Nevanlinna's Problem~\cite{rNicm}}.
Does $\underline E_m\ge 0$ imply that
the surface $X$ with the pullback complex structure
is parabolic?  Conversely, does $\overline E_m<0$ imply that it is hyperbolic?
\medskip

Teichm\"uller~\cite{oT38} constructed an example of a surface with
hyperbolic type, for which the mean excess is zero, thus
giving a negative answer to the first question.

We will shortly prove that the answer to the other question is negative
as well, by constructing a parabolic
surface $(\R^2,\psi)\in F_3$ with $E_m<0$.

In Section~\ref{S:nonpos}, we shall construct an example of a
non-positively curved, simply-connected, complete, parabolic
surface, whose curvature in any ball about a fixed basepoint is
less than a negative constant times the area of the ball.

\section{Counterexample}\label{S:construct}

P.~Doyle~\cite{pD84} proved that the surface $(X,\psi)$ is parabolic
if and only if a certain modification of the Speiser graph is recurrent.
(See~\cite{DS} and~\cite{pS94} for background on recurrence and
transience of infinite graphs.)
In the particular case where $k_f$ is bounded, the recurrence of the
Speiser graph itself is equivalent to $(X,\psi)$ being parabolic.
Though we will not really need this fact,
it is not too hard to see that in a Speiser graph satisfying $E_m<0$
the number of vertices in a ball grows exponentially with the radius.
Thus, we may begin searching for a counterexample by considering
recurrent graphs with exponential growth.
A very simple standard example of this sort is a tree constructed
as follows.
In an infinite $3$-regular tree $T_3$, let $v_0,v_1,\dots$ be
an infinite simple path.  Let $T$ be the set of vertices $u$
in $T_3$ such that $d(u,v_n)\le n$ for all sufficiently large $n$.
Note that there is a unique infinite simple path in $T$ starting from
any vertex $u$.  This implies that $T$ is recurrent.  It is
straightforward to check that the number of vertices
of $T$ in the ball $B(v_0,r)$ grows exponentially with $r$.

\begin{figure}
\tabskip=1em plus2em minus.5em
\halign to\hsize{\hfill#\hfill&\hfill#\hfill\cr
\includegraphics*[height=1.6in]{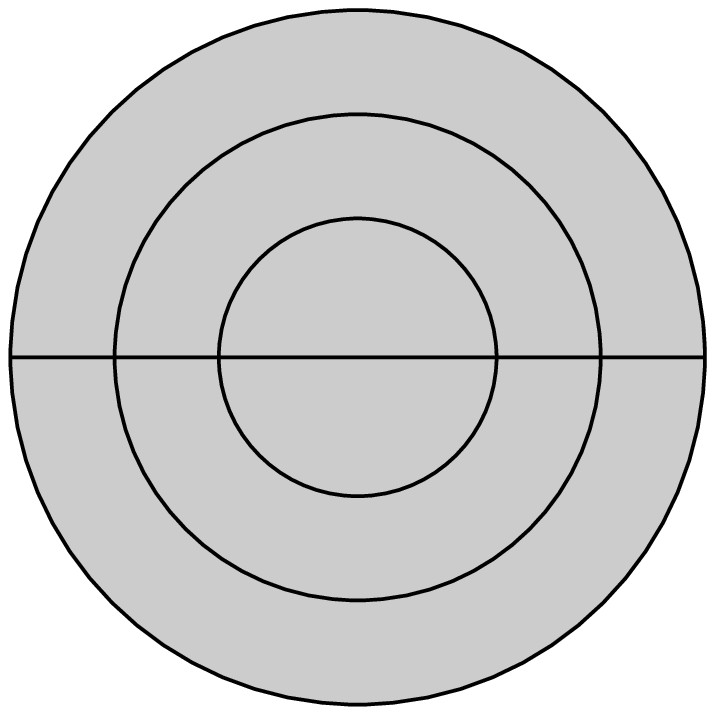}%
&%
\includegraphics*[height=1.6in]{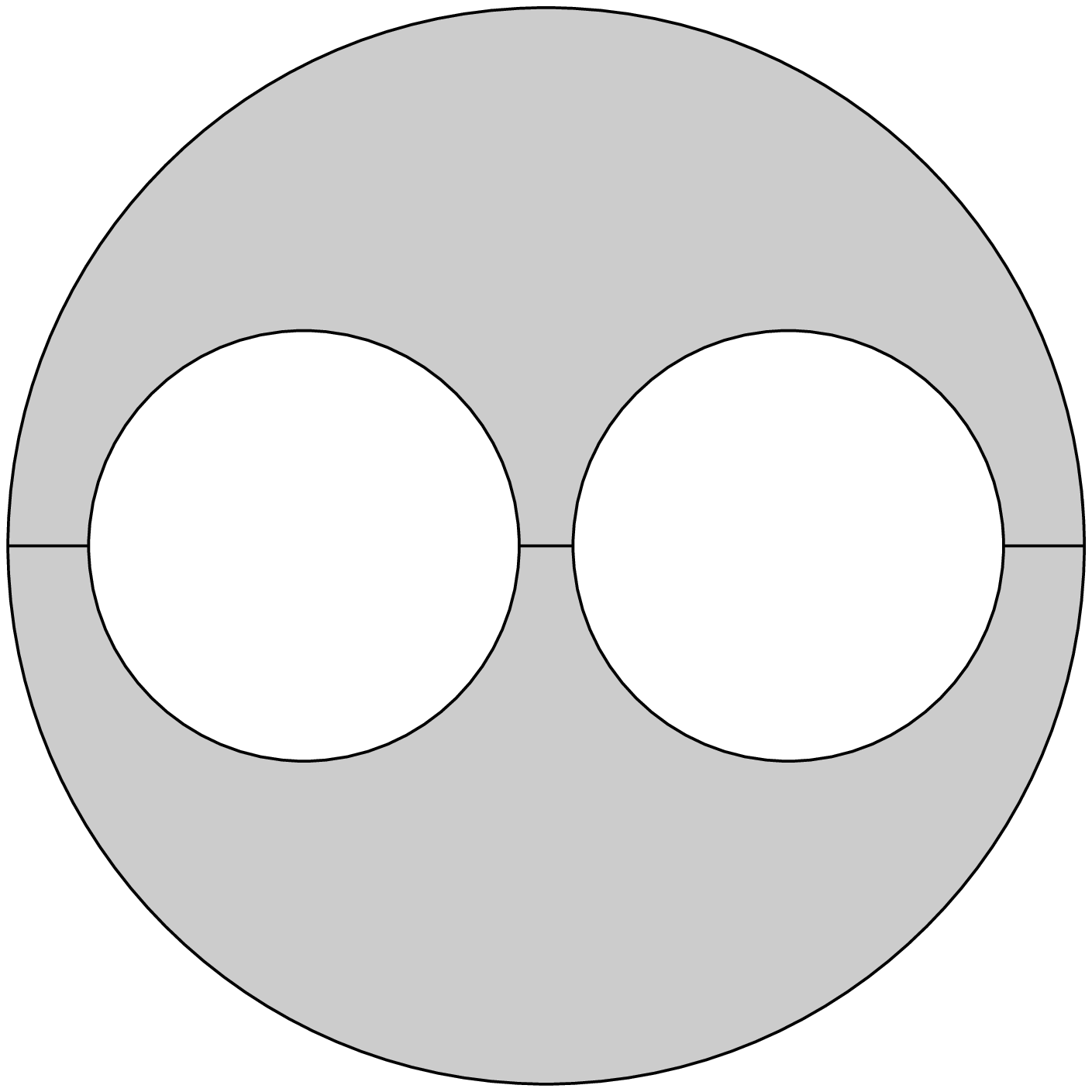}%
\cr
(a)&(b)\cr
}
\caption{\label{f.pantsandfinger}The surfaces $S(v)$.
(a) shows $S(v)$ for a leaf $v$, where $s=2$.
(b) shows $S(v)$ for a degree $3$ vertex.
}
\end{figure}

\begin{figure}
\SetLabels
(.29*.58)$S(v_1)$\\
(.4*.7)$S(v_2)$\\
(.65*.9)$S(v_3)$\\
\endSetLabels
\centerline{\AffixLabels{%
\includegraphics*[height=3in]{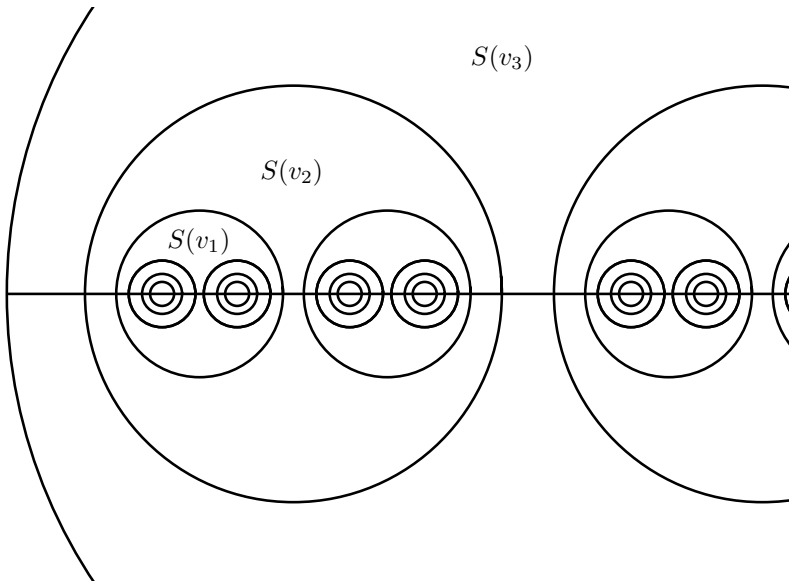}
}}
\caption{\label{f.Speiser}The Speiser graph with $s=2$.
The disk $S(v_0)$ is the left ``eye'' in $S(v_1)$.}
\end{figure}

Our Speiser graph counterexample
is a simple construction based on the tree $T$.
Fix a parameter $s\in\{1,2,\dots\}$, whose choice will be discussed
later.  To every leaf (degree one vertex) $v$ of $T$ associate a closed
disk $S(v)$ and on it draw the graph indicated in
Figure~\ref{f.pantsandfinger}.(a),
where the number of concentric circles, excluding
$\partial S(v)$, is $s$.
If $v$ is not a leaf, then it has degree $3$.
We then associate to it the graph indicated in
Figure~\ref{f.pantsandfinger}.(b),
drawn on a triply connected domain $S(v)$.  We combine these to form
the Speiser graph $\Gamma$ as indicated in figure~\ref{f.Speiser},
by pasting the outer boundary of the surface corresponding to each vertex
into the appropriate inner boundary component of its parent.
Here, the parent of $v$ is the vertex $v'$ such that
$d(v',v_n)=d(v,v_n)-1$ for all sufficiently large $n$.

Every vertex of $\Gamma$ has degree $4$ and every face
has $2$, $4$ or $6$ edges on its boundary.  Therefore,
$\Gamma$ is a Speiser graph.
Consequently, as discussed above, there is a surface
spread over the sphere $X=(\R^2,\psi)$ whose Speiser graph
is $\Gamma$.
It is immediate to verify that $\Gamma$ is recurrent,
for example, by the Nash-Williams criterion.
Doyle's Theorem~\cite{pD84} then implies that $X$ is parabolic.
Alternatively, one can arrive at the same conclusion
by noting that there is an infinite sequence of
disjoint isomorphic annuli on $(\R^2,\Gamma)$ separating
any fixed point from $\infty$, and applying extremal length.
(See~\cite{lA73,LV} for the basic properties of extremal
length.)

We now show that $\overline E_m<0$ for $\Gamma$.  Note that
the excess is positive only on vertices on the boundary
of $2$-gons, which arise from leaves in $T$.
On the other hand, every vertex of degree $3$ in $T$ gives
rise to vertices in $\Gamma$ with negative excess.  Take
as a basepoint for $\Gamma$ a vertex $w_0\in S(v_1)$ with negative excess,
say.  It is easy to see that there are constants,
$a>1,c>0$ such that the number $n^-_r$ of negative excess vertices
in the combinatorial ball $B(w_0,r)$ about $w_0$ satisfies
$c\,a^r\le n^-_r\le a^r/c$.

If $w$ is a vertex with positive excess, then there is a unique
vertex $\sigma(w)$ with negative excess closest to $v$;
in fact, if $w\in S(v)$, then $\sigma(w)$ is the closest
vertex to $w$ on $\partial S(v)$,
and the (combinatorial) distance from $w$ to $\sigma(w)$ is our parameter
$s$.  The map $w\mapsto \sigma(w)$ is clearly injective.
This implies that the number $n^+_r$ of positive excess
vertices in $B(w_0,r)$ satisfies
$n^+_{r+s}\le n^-_r$, $r\in\{0,1,2,\dots\}$.
By choosing $s$ sufficiently large, we may therefore
arrange to have the total excess in $B(w_0,r)$
to be less than $-\epsilon\, a^r$, for some $\epsilon>0$
and every $r\in\{0,1,2,\dots\}$.
It is clear that the number of vertices with zero excess in $B(w_0,r)$
is bounded by a constant (which may depend on $s$)
times $n^-_r$.  Hence, $\overline E_m<0$ for $\Gamma$.

By allowing $s$ to depend on the vertex in $T$, if necessary,
we may arrange to have $\underline E_m=\overline E_m$;
that is, $E_m$ exists, while maintaining $E_m<0$.
We have thus demonstrated that the resulting surface
is a counterexample in $F_4$ to the second implication in
Nevanlinna's problem.

\section{A non-positive curvature example}\label{S:nonpos}

We now construct an example of a simply connected,
complete, parabolic surface $Y$ of nowhere positive
curvature, with the property
\begin{equation}\label{E:Intcur}
\int_{D(a,r)}\text{curvature} <-\epsilon\,\area\bigl(D(a, r)\bigr),
\end{equation}
for some fixed $a\in Y$ and every $r>0$, where $D(a, r)$ denotes the
open disc centered at $a$ of radius $r$, and $\epsilon>0$
is some fixed constant.

Consider the surface $\C=\R^2$ with the metric $|dz|/y$ in
$P=\{z=x+iy:\ y\geq1\}$,
and $\exp{(1-y)}|dz|$ in $Q=\{y<1\}$. We denote this surface by $Y$.
Let $\beta$ denotes
the curve $\{y=1\}$ in $Y$, i.e., the common boundary of $P$ and $Q$.

Let $Q'$ denote the universal cover of $\{z\in\C:|z|>1\}$.
Note that $Q$ is isometric to $Q'$ via the map $z\mapsto \exp(iz+1)$.
Hence, the curvature is zero on $Q$, and the geodesic curvature
of $\partial Q$ is $-1$.  The geodesic curvature of $\partial P$ is $1$.
Consequently, $Y$ has no concentrated curvature on $\beta$.
The surface $Y$ is thus a ``surface of bounded curvature'',
also known as an Aleksandrov
surface (see~\cite{aA67}, \cite{yR93}).
The  curvature measure of $Y$ is absolutely
continuous with respect to area; the curvature of $Y$ is -1 (times area
measure) on $P$ and $0$ on $Q$.

The surface $Y$ is parabolic,
and the uniformizing map is the identity map onto $\R^2$ with the standard
metric.

We will now prove (\ref{E:Intcur}) with $a=i$.
Set $\beta_r=D(a,r)\cap\beta$. Note that the shortest path in $Y$ between
any two points on $\beta$ is contained in $P$, and is the arc
of a circle orthogonal to $\{y=0\}$. Using the Poincare disc model,
it is easy to see that there exists a constant $c>0$, such that
\begin{equation}\label{lbr}
c\,e^{r/2}\leq \length\beta_r\leq e^{r/2}/c,
\end{equation}
where the right inequality holds for all $r$, and the left for all sufficiently
large $r$.
By considering the intersection of $D(a,r)$ with the strip
$1<y<2$ it is clear that
\begin{equation}\label{E:Plb}
O(1)\area\bigl(P\cap D(a,r)\bigr)\ge \length\beta_r,
\end{equation}
for all sufficiently large $r$.

Consider some point $p\in Q$, and let
$p'$ be the point on $\beta$ closest to $p$.
It follows easily (for example, by using the
isometry of $Q$ and $Q'$) that if $q$ is any point in $\beta$,
then $d_Q(p,q)=d_Q(p,p')+d_Q(p',q)+O(1)$.
Consequently, if $d(p,a)\le r$,
then there is an $s\in[0,r]$ such that
$p'\in \beta_s$ and $d_Q(p,p')\le r-s+O(1)$.
Furthermore,
it is clear that the set of points $p$ in $Q$ such
that $p'\in \beta_s$ and $d_Q(p,p')\le t$ has
area $O(t^2+t)\,\length\beta_s$.  Consequently,
$$
\area\bigl(Q\cap D(a,r)\bigr)
\le O(1) \sum_{j=0}^{r} (j+1)^2\length\beta_{r-j}\,.
$$
Using~(\ref{lbr}), we have
\begin{equation}\label{E:Qub}
\area\bigl(Q\cap D(a,r)\bigr)\le O(1)\,\length\beta_r,
\end{equation}
for all sufficiently large $r$.

Now, combining~(\ref{E:Plb}) and (\ref{E:Qub}),
we obtain (\ref{E:Intcur}) for all sufficiently large $r$. It therefore
holds for all $r$.

\bigskip
\textbf{Acknowledgements}. {The authors are grateful to B.~Davis,
D.~Drasin, and especially to A.~Eremenko for helpful
discussions and their interest in this work.}

\enddocument